\documentclass[a4paper]{article}

\usepackage[english]{babel}
\usepackage[utf8x]{inputenc}
\usepackage[a4paper,top=3cm,bottom=2cm,left=3cm,right=3cm,marginparwidth=1.75cm]{geometry}
\usepackage{amsmath}
\usepackage{graphicx}
\usepackage[colorinlistoftodos]{todonotes}
\usepackage{listings}
\usepackage{color} 
\definecolor{mygreen}{RGB}{28,172,0} 
\definecolor{mylilas}{RGB}{170,55,241}
\usepackage{changepage}
\usepackage[small]{caption}
\usepackage{subcaption}
\usepackage[makeroom]{cancel}
\usepackage{amsthm}
\usepackage{amsfonts}
\usepackage{float}
\usepackage{framed}
\usepackage{cite}
\usepackage{algorithm}
\usepackage{pdfpages}
\usepackage{mathtools}
\usepackage{bm}
\usepackage{authblk}
\usepackage{multirow}
\usepackage{authblk}
\usepackage{tabularx}
\usepackage{pbox}
\usepackage{import}
\usepackage{tikz}
\usepackage{standalone}

\usepackage{todonotes}

\newcommand{\norm}[1]{\left\lVert#1\right\rVert}
\newcommand{\angstrom}{\text{\normalfont\AA}}

\newtheorem{remark}{Remark}

\begin{document}
	
	\title{Pseudo-Marginal Approximation to the Free Energy in a Micro-Macro Markov Chain Monte Carlo Method}
	
	\author[$\dag$]{Hannes Vandecasteele}
	\author[$\dag$]{Giovanni Samaey}
	\affil[$\dag$]{KU Leuven, Department of Computer Science, NUMA Section, Celestijnenlaan 200A box 2402, 3001 Leuven, Belgium}
	\date{\today}
	
	\maketitle
	
	\begin{abstract}
		We introduce a generalized micro-macro Markov chain Monte Carlo (mM-MCMC) method with pseudo-marginal approximation to the free energy, that is able to accelerate sampling of the microscopic Gibbs distributions when there is a time-scale separation between the macroscopic dynamics of a reaction coordinate and the remaining microscopic degrees of freedom. The mM-MCMC method attains this efficiency by iterating four steps: i) Propose a new value of the reaction coordinate; ii) Accept or reject the macroscopic sample; iii) Run a biased simulation that creates a microscopic molecular instance that lies close to the newly sampled macroscopic reaction coordinate value; iv) Microscopic accept/reject step for the new microscopic sample. In the present paper, we eliminate the main computational bottleneck of earlier versions of this method: the necessity to have an accurate approximation of the free energy. We show that introduction of a pseudo-marginal approximation significantly reduces the computational cost of the microscopic accept/reject step, while still providing unbiased samples. We illustrate the method's behaviour on several molecular systems with low-dimensional reaction coordinates.
	\end{abstract}
	
	\section{Introduction} \label{sec:introduction}
	The dynamics of large molecular systems, such as proteins, polymers and DNA, are determined by a potential energy surface $V(x)$ where $x \in \mathbb{R}^{3N}$ is the position vector of the $N$ atoms in the system. Molecules in a low energy state appear more frequently than molecules in a higher energy state. The likelihood of finding a molecule in a certain state $x$ is determined by the Gibbs distribution
	\begin{equation} \label{eq:gibbs}
		\mu(dx) = Z_V^{-1} \exp\left(-\beta V(x)\right) dx,
	\end{equation}
	with $Z_V$ the normalization constant,  $\beta = (k_B T)^{-1}$ the inverse temperature and $k_B$ the Boltzmann constant in appropriate units. The main challenge in computational chemistry is to find the local minima of $V(x)$, by means of sampling~\eqref{eq:gibbs}.
	
	Sampling~\eqref{eq:gibbs} is especially hard when there is metastability or a time-scale separation between the microscopic, high-dimensional atomistic positions, and the macroscopic degrees of freedom. These macroscopic degrees of freedom largely determine the structure of the molecule. The macroscopic variables are usually represented by a reaction coordinate
	\begin{align}
		\xi &: \mathbb{R}^{3N} \to \mathbb{R}^m \\
		&x \mapsto z = \xi(x).
	\end{align}
	When the reaction coordinate represents the true macroscopic dynamics, the potential energy can be decomposed as
	\begin{equation} \label{eq:decomposition}
		V(x) = \frac{1}{\varepsilon} V_f(x) +  A(\xi(x)), 
	\end{equation}
	with $1/\varepsilon$ the time-scale separation of the system. Here, $V_ f(x)$ describes the dynamics of the microscopic variables, and $A(\xi(x))$ that of the reaction coordinate. Both are assumed independent from $\varepsilon$. The function $A \circ \xi(x)$ is the free energy of the reaction coordinate~\cite{stoltz2010free} defined by
	\begin{equation}
		A(z) = -\beta^{-1} \ln \int_{\Sigma(z)} \mu(x) \ \delta_{\xi(x) - z}(dx),
	\end{equation}
	with $\Sigma(z)  = \{ x \in \mathbb{R}^{3N} \ | \ \xi(x) = z  \}$ the level set of the reaction coordinate at $z$. It is opportune at this time to define the time-invariant distribution of reaction coordinate values
	\begin{equation}
		\mu_0(dz) = Z_A^{-1} \exp\left(-\beta A(z)\right) dz,
	\end{equation}
	with normalization constant $Z_A$.

	We have recently developed a micro-macro Markov chain Monte Carlo method (mM-MCMC) that enhances sampling of~\eqref{eq:gibbs}~\cite{vandecasteele2022} in case of metastability and/or time-scale separation. The mM-MCMC method improves exploration of the microscopic state space $\mathbb{R}^{3N}$ by first generating samples on the macroscopic level, and then going back to the microscopic state space by means of reconstruction. The latter is the hardest part of the algorithm. We have developed two techniques: direct and indirect reconstruction. Only indirect reconstruction is used in this manuscript. Putting these together, the mM-MCMC method consists of four steps; i) \textbf{Macroscopic Proposal}: Generate a new reaction coordinate value using a proposal distribution; ii) \textbf{Macroscopic Acceptance}: Accept or reject this reaction coordinate value using an approximate macroscopic distribution; iii) \textbf{Indirect Reconstruction}: Run a biased simulation that pulls the microscopic sample towards the level set of the new reaction coordinate value; iv) \textbf{Microscopic Acceptance}: Accept or reject the newly computed microscopic sample according to the Gibbs distribution~\eqref{eq:gibbs}.
	
	We showed in previous work that this method converges to the Gibbs distribution~\cite{vandecasteele2022}, and that it can obtain significant efficiency gains over traditional microscopic samplers~\cite{vandecasteele2023}. This efficiency gain comes at a cost, however. Knowing the free energy of the reaction coordinate is crucial to obtain unbiased samples of~\eqref{eq:gibbs}. As we will cover in section~\ref{sec:mMMCMC}, the free energy is required to evaluate the microscopic acceptance probability. Any bias on the free energy results in an error on the Gibbs distribution, and we cannot measure how large the resulting bias can be. Solving the problem of the free energy makes mM-MCMC an efficient alternative for multiscale sampling applications.
	
	We give a short summary of approximation schemes to the free energy in the literature. The adaptive biasing force/potential methods are a general class of schemes~\cite{dickson2010free,comer2015adaptive,lelievre2008long} that compute a local approximation to the respective gradient of free energy, or the free energy itself. By subtracting the force/potential from the microscopic dynamics creates a flat energy landscape in the direction of the reaction coordinate. Another class of methods simulates a non-equilibrium dynamics between reaction coordinate value and computes the free energy difference by the Jarzynski theorem~\cite{lelievre2007computation,lelievre2012langevin}. Other methods use thermodynamic integration, usually constrained to the manifold of given reaction coordinate value. Sampling is typically done using a Rattle scheme~\cite{andersen1983rattle} or generalized Hamiltonian Monte Carlo (GHMC)~\cite{lelievre2012langevin}.
	
	In this manuscript, we use the pseudo-marginal approach~\cite{andrieu2009pseudo,filippone2014pseudo,alenlov2016pseudo,sherlock2015efficiency} to approximate the free energy at every iteration of the mM-MCMC algorithm. This approximation generates microscopic samples from the higher dimensional distribution and marginalises them, to compute the free energy using importance sampling. Key is that this method is unbiased and can be run on the fly at every mM-MCMC step. The only downside is an added variance because of  the importance sampling.
	
	\paragraph{Outline}
	This manuscript follows the subsequent structure. We explain the mM-MCMC method with indirect reconstruction in section~\ref{sec:mMMCMC}, together with a discussion of its invariant distribution. We also discuss the problem of the free energy in this context. Then in section~\ref{sec:three_fec}, we outline three free energy computation methods from the literature, explain them in some detail, and discuss their use in the context of mM-MCMC. Afterwards, section~\ref{sec:pseudomarginal} outlines the pseudo-marginal method for sampling marginal distributions, and we adapt this scheme to the mM-MCMC reconstruction step in section~\ref{sec:pmreconstruction}. Finally, we test the pseudo-marginal mM-MCMC method on butane and N-alkane in section~\ref{sec:numerics}. We conclude the manuscript in section~\ref{sec:conclusion}.

	\section{The Micro-Macro MCMC method and Invariant Distribution} \label{sec:mMMCMC}
	This section serves to outline the mM-MCMC method, together with the free energy problem. We explain the steps of the mM-MCMC method in section~\ref{subsec:mMMCMC}. We then define its time-invariant probability distribution in section~\ref{subsec:invariant_mu}. Afterwards, we show how the mM-MCMC method depends on the free energy of the reaction coordinate in section~\ref{subsec:exactness} and how it affects the performance of the mM-MCMC method.
	
	\subsection{Micro-Macro Markov Chain Monte Carlo} \label{subsec:mMMCMC}
	The mM-MCMC method aims to accelerate and improve sampling of the Gibbs distribution~\eqref{eq:gibbs} in case of metastability or multiple time scales. We introduced the mM-MCMC method in~\cite{vandecasteele2022} and compared it to other schemes in~\cite{vandecasteele2023}. The first paper contains two variants of the mM-MCMC scheme: with direct and indirect reconstruction. We will only use the latter in this manuscript, as the former is not always applicable.
	
	Given a microscopic sample $x_n \in \mathbb{R}^{3N}$ and a reaction coordinate value $z_n \in \mathbb{R}^m$ that doesn't necessarily equal to $\xi(x_n)$. The mM-MCMC method with indirect reconstruction computes the next sample $(x_{n+1}, z_{n+1})$ in four steps.
	
	\subsubsection{Macroscopic Proposal}
	Generate a new reaction coordinate value $z'$ with the macroscopic proposal distribution $q_0(z_n, z')$. These proposals are usually based on a time-stepping scheme such as Brownian Motion, gradient descent ore more complex methods such as splitting schemes and Hamiltonian Monte Carlo.
	
	Since the macroscopic system is of low dimensionality and not stiff, we are able to use large time steps of size $\Delta t$, which also helps overcome metastability.
	
	\subsubsection{Macroscopic Acceptance}
	Before proceeding to the microscopic level, we must make sure that $z'$ samples an approximate macroscopic distribution $\bar{\mu}_0$. We accept $z'$ with probability
	\begin{equation}\label{eq:macro_acc}
		\alpha_{CG} (z_n, z') =  \min\left\{1,  \frac{\bar{\mu}_0(z')  q_0(z', z_n)}{\bar{\mu}_0(z_n)  q_0(z_n, z')} \right \}.
	\end{equation}
	On acceptance, we proceed to reconstruction, else $(x_{n+1}, z_{n+1}) = (x_n, z_n)$. Note that we cannot use the free energy as our approximate macroscopic distribution because direct evaluations are very expensive.
	
	\subsubsection{Indirect Reconstruction}
	With $z'$, we create a microscopic sample $x'$ such that $\xi(x')  \approx z'$.  Starting from the previous microscopic sample $x_n = x_n^0$, we run $K$ steps of the biased simulation
	\begin{equation} \label{eq:biased_simulation}
		x_n^{k+1} = x_n^k  - \delta t \nabla V(x_n^k) - \delta t \lambda (\xi(x_n^k) - z')  \nabla \xi(x_n^k) + \sqrt{2 \delta t  \beta^{-1}} \ \eta_k, \ \ k = 0, \dots, K-1, 
	\end{equation}
	and set $x' = x_n^K$. Here, $\delta t$ is the reconstruction time step and $\eta_k$ are independent standard normally distributed numbers. The parameter $\lambda$ biases $x_n^k$ near the level set $\Sigma(z')$. If $\lambda$ is large, it also makes the slow variable $\xi$ fast, and on the same time scale as the other variables. Therefore, choosing $\lambda$ on par with the fastest mode of the system creates good mixing of~\eqref{eq:biased_simulation}~\cite{vandecasteele2023}.
	
	The biased process~\eqref{eq:biased_simulation} has an invariant distribution
	\begin{equation} \label{eq:nu_laambda}
		\nu_\lambda(x  | z') = \mu_\lambda(z')^{-1} \exp\left(-\beta V(x)\right) \exp\left(-\beta \lambda \frac{\norm{\xi(x) - z'}^2}{2} \right),
	\end{equation}
	where the normalization constants is
	\begin{equation} \label{eq:N_lambda}
		\mu_\lambda(z')  = Z_V \int_{\mathbb{R}^n} \left(\frac{\lambda  \beta}{2 \pi}\right)^{n/2} \exp\left(-\beta \lambda \frac{\norm{u - z'}^2}{2} \right) d\mu_0(u).
	\end{equation}

	\subsubsection{Microscopic Acceptance}
	Finally, to guarantee that $x'$ samples the Gibbs distribution, we accept or reject $x'$ with probability
	\begin{equation} \label{eq:micro_acc}
		\alpha_f(z_n, x_n, z', x') = \min\left\{1, \frac{\mu(x') \bar{\mu}_0(z_n) \nu_\lambda(x_n|z_n)}{\mu(x_n) \bar{\mu}_0(z') \nu_\lambda(x'|z')} \right\} = \min \left\{1, \frac{\mu_\lambda(z') \bar{\mu}_0(z_n)}{\mu_\lambda(z_n) \bar{\mu}_0(z')} \right\}.
	\end{equation}
	On acceptance, $(x_{n+1}, z_{n+1}) = (x', z')$, otherwise  $(x_{n+1}, z_{n+1}) = (x_n, z_n)$.
	
	\subsection{Invariant Distribution} \label{subsec:invariant_mu}
	Due to the nature of the indirect reconstruction scheme, there is no direct coupling between the microscopic sample $x'$ and the reaction coordinate value $z'$, i.e., $\xi(x') \neq z'$. The same is true for $x_n$ and $z_n$. As a consequence, the mM-MCMC scheme does not sample the Gibbs distribution only, but an extended distribution on the larger state space of $\mathbb{R}^m \times \mathbb{R}^{3N} $. This distribution reads
	\begin{equation} \label{eq:mu_extended}
		\mu_{\text{ext}}(z,x) = \left(\frac{\lambda \beta}{2 \pi}\right)^{n/2}  \exp\left(-\beta \lambda \frac{\norm{\xi(x) - z}^2}{2} \right) \mu(x),
	\end{equation}
	where marginal in $x$ of $\mu_{\text{ext}}(x, z)$ is the Gibbs distribution. We already showed~\cite{vandecasteele2022} that the mM-MCMC scheme indeed samples the above distribution and that the method is ergodic.
	
	Moreover, one can also write~\eqref{eq:mu_extended} as
	\begin{equation}
		\mu_{\text{ext}}(z,x) = \mu_\lambda(z) \nu_\lambda(x | z),
	\end{equation}
	with $\mu_\lambda$ the marginal distribution in $z$, and $\nu_\lambda(x|z)$ the conditional distribution of $x$, given $z$. We will need this expression later.
	
	\subsection{Exactness and Free Energy} \label{subsec:exactness}
	The expression of the microscopic acceptance probability~\eqref{eq:micro_acc} requires an evaluation of $\mu_\lambda$ after every reconstruction step. Recalling the definition~\eqref{eq:N_lambda}, $\mu_\lambda$ is given by an integral over the marginal invariant distribution of the reaction coordinate. We can never evaluate $\mu_\lambda$ exactly, nor can we pre-compute it because of the bias it would induce on the whole mM-MCMC sampling result. The only way to compute the microscopic acceptance probability is to rely on an on-the-fly estimate $\tilde{\mu}_\lambda$ of $\mu_\lambda$ that maintains exactness of the mM-MCMC algorithm.
	
	At this time, it is opportune to introduce the partition function $Q_\lambda(z)$ of $\mu_\lambda$
	\begin{equation} \label{eq:Q_lambda}
		Q_\lambda(z) = -\frac{1}{\beta} \ln \mu_\lambda(z),
	\end{equation}
	so that the microscopic acceptance probability can be written as
	\begin{equation} \label{eq:micro_acc_variance}
		\alpha_f(z_n, x_n, z', x') = \min \left\{1, e^{-\beta \left(Q_\lambda(z') - Q_\lambda(z_n)\right)} \frac{ \bar{\mu}_0(z_n)}{ \bar{\mu}_0(z')} \right\}.
	\end{equation}
	From the definition~\eqref{eq:Q_lambda} , $Q_\lambda$ acts as the `free energy' of the marginal distribution $\mu_\lambda$ of a new reaction coordinate on the extended state space $\mathbb{R}^{3N} \times \mathbb{R}^m$ given by
	\begin{align} \label{eq:new_rc}
		\xi_{\text{ext}} : &  \mathbb{R}^m \times \mathbb{R}^{3N}   \to \mathbb{R}^m\\
		&(z,x) \mapsto z = \xi_{\text{ext}}(z,x).
	\end{align}
	It is this free energy we need to compute to evaluate the microscopic acceptance probability.
	
	In the remainder of this manuscript, whenever we mention the `free energy' or the `marginal distribution', we refer to $Q_\lambda$ and $\mu_\lambda$ respectively.
	
	\section{An Overview of Free Energy Computations} \label{sec:three_fec}
	Free energy computations (FEC) are central in the multiscale molecular dynamics literature. We give a couple of key references~\cite{stoltz2010free,homeyer2012free,chipot2007free}. We do not need to approximate the free energy up to double precision in most applications. Indeed, only a rough free energy difference estimate 
	\begin{equation*}
		\tilde{Q}_\lambda(z') - \tilde{Q}_\lambda(z_n)
	\end{equation*}
	is enough to evaluate~\eqref{eq:micro_acc_variance}, as long the approximation is exact in expectation. In what follows we give a summary of the most important free energy computation schemes from the literature. We present each scheme in some detail, and discuss its applicability for evaluating $Q_\lambda$ in the context of the microscopic acceptance probability in the mM-MCMC algorithm. We end this section by presenting a list of properties that any FEC scheme should satisfy to estimate $Q_\lambda$ in the context of~\eqref{eq:micro_acc_variance}.
	
	\subsection{Three Free Energy Computation Schemes from the Literature}
	The three free energy computation schemes that we will present in this section are the biased potential / force (ABP/ ABF) methods, Thermodynamic Integration (TI), and a non-equilibrium stochastic scheme based on the famous Jarzynski relation~\cite{jarzynski1997nonequilibrium}. We discuss where each of these methods work, and where they do not, and bring then in relation with~\eqref{eq:micro_acc_variance}.
	
	\subsubsection{Adaptive Biasing Potential  Force}
	The idea behind the Adaptive Biasing Potential (ABP) method is to build up knowledge of the free energy $Q_\lambda(z)$ over time, and then to subtract it from the microscopic process to `flatten' the microscopic dynamics in the direction of the reaction coordinate. Suppose we are interested in the value of $\mu_\lambda$ or $Q_\lambda$ at the value $z^*$. The microscopic process in question is typically the overdamped Langevin dynamics that keeps $\mu_{\text{ext}}$ invariant
	\begin{equation} \label{eq:overdamped_langevin}
		dq_t = \nabla \ln \left(\mu_{\text{ext}} (q_t) \right)dt + \sqrt{2\beta^{-1}} dW(t), 
	\end{equation}
	with $q_t = (z_t, x_t)$. 
	
	The Adaptive Biased Potential method then approximates the true free energy $Q_\lambda$ at $z^*$ by a mollification $Q_\lambda^{\varepsilon}(z^*)$
	\begin{equation} \label{eq:ABP_mollification}
		\exp\left(-\beta  Q_\lambda^{\varepsilon}(z^*, t) \right)  = Z_t^{-1} \left(1 + \int_0^t \delta_\varepsilon (z_s - z^*) ds\right),
	\end{equation}
	where $\delta_\varepsilon$ is a Gaussian approximation to the delta function with variance $\varepsilon^2$, and $Z_t$ is the normalization constant so that $Z_t^{-1} \int_{\mathbb{R}^n} \exp\left(-\beta Q_\lambda^{\varepsilon}(z, t) \right) dz = 1$. Equation~\eqref{eq:ABP_mollification} aggregates the macroscopic samples $z_s$ and weighs them according to their distance with respect to $z^*$. Close sample have a higher weighing than farther samples.
	
	The Adaptive Biased Force (ABF) method takes a different approach by computing the gradient $\nabla_{z^*} Q_\lambda^\varepsilon (z^*)$ instead,
	\begin{equation}
		\nabla_{z^*} Q_\lambda^{\varepsilon}(z^*, t) = -\beta^{-1} \frac{\int_0^t\nabla_{z^*} \delta_\varepsilon(z_s - z^*) ds}{1 + \int_0^t \delta_\varepsilon(z_s - z^*) ds}.
	\end{equation}
	We refer for more details to~\cite{lelievre2008long,comer2015adaptive,dickson2010free,darve2008adaptive} about the theory and implementation.
	
	We have have two remarks concerning the ABP/ABF framework. First, these methods work best when one wants to approximate the free energy at many points at the same time. The Langevin dynamics~\eqref{eq:overdamped_langevin} is independent of $z^*$, and the update procedure~\eqref{eq:ABP_mollification} can be implemented in parallel.
	
	Second, both methods rely on the overdamped Langevin dynamics to generate microscopic samples that are generated according to~\eqref{eq:mu_extended}. If we were to use the ABP method to calculate the free energy difference in the mM-MCMC algorithm, the indirect reconstruction and overdamped Langevin simulations would need to run simultaneously, adding overhead. We do not currently see use for the ABP/ABF methods for this reason. 
	
	\subsubsection{Thermodynamic Integration}
	The method of thermodynamic integration is based on two observations~\cite{chipot2007free,hummer2001fast,stoltz2010free}. The first observation is that the free energy difference between two states $z$ and $z'$ can be written in terms of the free energy derivative $d Q_\lambda / d \xi_{\text{ext}}$
	\begin{equation} \label{eq:TI}
		Q_\lambda(z') - Q_\lambda(z) = \int_{z}^{z'} \frac{dQ_\lambda}{ d\xi_{\text{ext}}} d \xi_{\text{ext}},
	\end{equation}
	and that the free energy derivative, in turn, can be computed using a constrained microscopic statistical average
	\begin{equation} \label{eq:TI_derivative}
		\frac{d Q_\lambda}{ d \xi_{\text{ext}}}(z) = \frac{\int_{\mathbb{R}^{3N}} \frac{\partial V_{\text{ext}}(z,x)}{\partial z} \mu_{\text{ext}} (z,x) \ dx }{\int_{\mathbb{R}^{3N}} \mu_{\text{ext}} (z,x) \ dx } = \int_{\mathbb{R}^{3N}} \frac{\partial V_{\text{ext}}(z,x)}{\partial z} \mu_{\text{ext}} (z,x) \ dx. 
	\end{equation}
	This is the second observation. In notation above, $V_{\text{ext}}$ is the potential energy of $\mu_{\text{ext}}$. In a practical implementation, one would choose $M$ collocation points $z_m$ between $z$ and $z'$, approximate the free energy derivative $\frac{d Q_\lambda}{ d \xi_{\text{ext}}}(z_m) $ at each such point, and then integrate~\eqref{eq:TI} using a quadrature rule.
	
	We see two reasons why thermodynamic integration cannot be used efficiently to compute the free energy difference in the mM-MCMC algorithm. The first reason is about computation cost. For an accurate estimate of the free energy derivative, we need many particles over a large part of the microscopic state space $\mathbb{R}^{3N}$. This is too high a cost for TI to be used after every reconstruction step. The second reason is about accuracy. Even if one can compute the free energy derivative to any specified precision, there will always be a deterministic error on~\eqref{eq:TI}. Error control typically depends on a higher order quadrature rule, making TI too expensive
	
	\subsubsection{Non-Equilibrium Integration}
	The final FEC method we present is a non-equilibrium constraint dynamics together with the renowned Jarzynski formula~\cite{jarzynski1997nonequilibrium}. The Jarzynski formula relates the free energy difference between two states $z$ and $z'$ to the irreversible work done by an ensemble of constraint paths joining $z$ with $z'$.
	
	For a practical implementation, we first need a monotone differentiable schedule of reaction coordinate values $z(t)$ so that $z(0) = z$ and $z(T) = z'$ for some end time $T \geq 0$. The constrained microscopic Langevin dynamics then reads
	\begin{equation} \label{eq:constrained_langevin}
		\begin{cases}
			q_0 &\sim \nu_\lambda(\cdot | z) \\
			dq_t &= -\nabla V_{\text{ext}} (q_t) dt + \sqrt{2\beta^{-1}} dW_t + \sum_{\alpha=1}^m \nabla \xi_{\text{ext}, \alpha}(q_t) d\lambda _{\alpha, t} \\
			z(t) &= 	\xi_{\text{ext}}(q_t).
		\end{cases}
	\end{equation}
	In the equation above, $\lambda_{\alpha, t}$ are Lagrange multipliers associated to the $m$ constraints in the third equation. The mechanical work done by the constrained dynamics along one path $q_t, t = 0, \dots, T$ can be computed by
	\begin{equation*}
		\mathcal{W}_{0, T}(q_t) = \sum_{\alpha=1}^m \int_0^T f_\alpha(q_s) \dot{z}_\alpha(s) ds.
	\end{equation*}
	The function $f_\alpha$ is the local mean force of the reaction coordinate, see~\cite{lelievre2007computation} for more details. Finally, the Jarzynski formula relates the average work by the constrained simulation to the free energy difference by
	\begin{equation}
		e^{-\beta \left(Q_\lambda(z') - Q_\lambda(z) \right)} = \mathbb{E}\left[e^{-\beta\mathcal{W}_{0, T}(q_t)}\right].
	\end{equation}
	We refer to~\cite{stoltz2010free,lelievre2012langevin,lelievre2007computation,rousset2006equilibrium} for more details regarding discretizations of this scheme. 
	
	There are many similarities between the indirect reconstruction scheme and the constrained process. Both `pull' the previous microscopic sample $x_n$ (or $q_0$) to the new reaction coordinate value, although indirect reconstruction does not enforce this constraint. The downside of the constrained / non-equilibrium process is many paths are needed to obtain a precise estimate of the free energy difference. We have also found the non-equilibrium scheme to be slow due to the many projections during the simulation of~\eqref{eq:constrained_langevin}.
	
	\subsection{Discussion}
	We end this section by going a list of necessary requirements for FEC methods regarding their use for computing the free energy in the mM-MCMC algorithm. These requirements are the following:
	\begin{itemize}
		\item Speed: The FEC method should be fast and have little overhead because we compute the free energy after every reconstruction step.
		
		\item Work well with Indirect Reconstruction: Indirect reconstruction is expensive itself, so the FEC method ideally use the same reconstructed samples to generate an approximation.
		
		\item Exact in Expectation: The FEC estimator must be exact in estimation, or any deterministic error must be controllable in term of the parameters of the FEC method.
	\end{itemize}
	We present a different FEC method that fulfils each of these criteria in the next section.
	
	\section{The Pseudo-Marginal Approximation} \label{sec:pseudomarginal}
	Computing the marginal distribution $\mu_\lambda$ from a joint distribution $\mu_{\text{ext}}$  is a central problem in many applications~\cite{beaumont2003estimation}. A relatively easy and tractable method to sample from $\mu_\lambda(z)$ would be to use a Markov chain Monte Carlo (MCMC) method that generates samples $\{(z_i,x_i)\}$ from the joint distribution $\mu_{\text{ext}}(z, x)$, and then single out the macroscopic particles $\{z_i\}$. However, it is well established that because of stiffness in~\eqref{eq:decomposition}, such MCMC methods can result in strongly correlated samples, an undesirable property. On the other hand, if $\mu_\lambda(z)$ were analytically known or was easy to compute, one might simply generate samples from $\mu_\lambda(z)$ in a Markov chain Monte Carlo fashion. Indeed, if the Markov chain is at position $z$, the MCMC method generates a new sample $z′$ using some proposal distribution $q(z,z′)$. We accept this proposal with probability
	\begin{equation} \label{eq:alpha}
		\alpha(z, z') = \min \left\{1, \frac{\mu_\lambda(z') q(z', z)}{\mu_\lambda(z) q(z, z')} \right\}
	\end{equation}
	which could be evaluated analytically.

	The pseudo-marginal method has been designed to combine the benefits of both approaches~\cite{andrieu2009pseudo,alenlov2016pseudo,sherlock2015efficiency,andrieu2015convergence}:
	possible computational and statistical efficiency gains by directly sampling from $\mu_\lambda(z)$ and ease of implementation by using auxiliary microscopic variables that sample from $\mu_{\text{ext}}(z,x)$. A particularly natural approach for estimating the intractable marginal density $\mu_\lambda(z)$ is importance sampling. That is, for some integer $K \geq 1$ and some conditional distribution $q_z$, we consider the estimates
	\begin{align} \label{eq:IS}
		\tilde{\mu}_\lambda(z) &= \frac{1}{K} \sum_{k=1}^K \frac{\mu_{\text{ext}}(z, x_k)}{q_z(x_k)},  \ x_k|z \sim q_z, \text{iid}, \\
		\tilde{\mu}_\lambda(z') &= \frac{1}{K} \sum_{k=1}^K \frac{\mu_{\text{ext}}(z', x'_k)}{q_z'(x'_k)},  \ x'_k|z' \sim q_z', \text{iid}
	\end{align}
	These estimates are unbiased for $\mu_\lambda$ as long as all microscopic variables are independent. Let us denote the microscopic samples of~\eqref{eq:IS} by $\mathcal{X} = \{x_1, . . . , x_k\}$ and $\mathcal{X}' = \{x'_1,..., x'_K\}$. These microscopic variables then have densities $q_z(\mathcal{X})$ and $q_{z'}(\mathcal{X}')$, respectively. 
	
	With the above estimates for the marginal distribution $\mu_\lambda(z)$ we can complete the Markov chain that we introduced at the top of the section. Specifically, if we plug $\tilde{\mu}_\lambda(z) $ and $\tilde{\mu}_\lambda(z')$ into acceptance probability~\eqref{eq:alpha}, we obtain
	\begin{equation*}
		\tilde{\alpha}(z, z') = \min\left\{1, \frac{\left[\frac{1}{K} \sum_{k=1}^K \frac{\mu(z', x'_k)}{q_z'(x'_k)} \right] q(z', z)}{\left[\frac{1}{K} \sum_{k=1}^K \frac{\mu(z, x_k)}{q_z(x_k)} \right] q(z, z')}\right\}.
	\end{equation*}
	However, there is no guarantee that we are still sampling from the marginal distribution $\mu_\lambda(z)$. A remarkable property~\cite{beaumont2003estimation} is that we can rewrite the above (approximate) acceptance rate as
	\begin{equation} \label{eq:alpha_pm}
		\tilde{\alpha}(z, z') = \min\left\{1, \frac{\left[ \frac{1}{K} \sum_{k=1}^K \mu_{\text{ext}}(z', x_k') \Pi_{l=1, l\neq k}^K q_{z'}(x'_l)\right] q(z', z) q_z^K(\mathcal{X})}{\left[ \frac{1}{K} \sum_{k=1}^K \mu_{\text{ext}}(z, x_k) \Pi_{l=1, l\neq k}^K q_{z}(x_l)\right] q(z, z') q_z'^K(\mathcal{X}')}\right\}.
	\end{equation}
	This rearrangement of the (estimated) acceptance rate suggests that the pseudo-marginal method is a Markov chain on the larger state space $(z,\mathcal{X})$, not just in $z$. Indeed, the pseudo-marginal method constructs a Markov chain with transition probability $q(z, z' ) \ q_{z′} (\mathcal{X} )$ and with invariant measure what is in between the square brackets. This formula implies that the marginal distribution in $z$ of
	\begin{equation*}
		\tilde{\mu}(z, \mathcal{X}) =  \frac{1}{K} \sum_{k=1}^K \mu_{\text{ext}}(z, x_k) \ \Pi_{l=1, l\neq k}^K q_{z}(x_l),
	\end{equation*}
	is still $\mu_\lambda(z)$. That is, one is still guaranteed to sample from the correct marginal distribution.
	
	It is possible to further optimize the pseudo-marginal method by noting that we must only compute $\tilde{\mu}_\lambda(z')$ once. Indeed, on acceptance of $z'$ a new macroscopic sample is generated and the exact same term now appears in the denominator of~\eqref{eq:alpha_pm}. This adaptation of the pseudo-marginal scheme is called the Grouped Independence Metropolis-Hastings method (GIMH)~\cite{beaumont2003estimation,andrieu2009pseudo}, and is used in all our implementations of the method.

	\section{Adapting the Pseudo-Marginal Approach to Reconstruction} \label{sec:pmreconstruction}
	Having established the pseudo-marginal method as an algorithm for sampling from and evaluating the marginal distribution, we now turn our attention to approximating the normalization constant $\mu_{\lambda}$ inside the microscopic acceptance probability~\eqref{eq:micro_acc}. In section~\ref{subsec:histogram}, we explain the main changes to the pseudo-marginal method to tailor it specifically to the mM-MCMC method. Then, in section~\ref{subsec:tensor}, we present a faster alternative for practical implementations. Afterwards we illustrate the full mM-MCMC method with the pseudo-marginal approximation in section~\ref{subsec:algorithm}.
	
	\subsection{A Histogram Approximation to $q_z$} \label{subsec:histogram}
	During indirect reconstruction, a sequence $\mathcal{X}_K = \{x_k\}_{k=1}^K$ of microscopic samples is constructed from reconstruction distribution $\nu_{\lambda}$, at least asymptotically. It therefore makes sense to choose $q_z = \nu_{\lambda}(\cdot|z)$ as our importance distribution in the pseudo-marginal method. The reconstructed microscopic samples are then used to estimate the free energy $\mu_{\lambda}$ using~\eqref{eq:IS}. There is thus an intimate connection between indirect reconstruction and pseudo-marginalization.
	
	There is one more problem, however. The (filtered) free energy $\mu_\lambda$ is the normalization constant of the reconstruction distribution $\nu_\lambda$. That is, we need to evaluate $\mu_\lambda$ to write an estimator for $\mu_\lambda$, which is nonsensical. We can overcome this issue by not using $\nu_{\lambda}$ as $q_z$, but an approximation to it that has already been normalized. A histogram is the ideal tool for this.
	
	\subsubsection{Constructing the Histogram}
	To create a histogram $H_z[\mathcal{X}_K]$ that approximates $\nu_\lambda$, we discretize $\mathbb{R}^{3N}$ into equally sized square blocks or `bins'. Suppose $h$ is the bin length, so that $h^{3N}$ is the total volume of one bin. Before actually constructing the histogram, let us introduce some notation. A particle $q \in \mathbb{R}^{3N}$ falls into bin $\left(l_1, l_2, \dots, l_{3N-1}, l_{3N}\right) \in \mathbb{Z}^{3N}$ if
	\begin{equation} \label{eq:bin_condition}
		l_i h \leq q^{(i)} < (l_i+1)h,
	\end{equation}
	for all $i$ from $1$ to $3N$. We use the notation $q^{(i)}$ to denote the $i$-th component of $q$. In equation~\eqref{eq:bin_condition}, we implicitly assume that the origin is a corner of the `central' bin.
	
	Let $K_{\left(l_1, l_2, \dots, l_{3N}\right)}$ be the number of reconstructed samples that fall into bin $\left(l_1, l_2, \dots, l_{3N}\right)$, where 
	\begin{equation*}
		\sum_{l_1, l_2, \dots, l_{3N} \in \mathbb{Z}} K_{\left(l_1, l_2, \dots, l_{3N}\right)} = K.
	\end{equation*}
	We then construct the histogram $H_z[\mathcal{X}_K](x)$ as
	\begin{equation} \label{eq:3N_histogram}
		H_z[\mathcal{X}_K](x) = \sum_{l_1, l_2, \dots, l_{3N} \in \mathbb{Z}} \mathcal{I}_{\left(l_1, l_2, \dots, l_{3N}\right)}(x) \frac{K_{\left(l_1, l_2, \dots, l_{3N}\right)}}{K},
	\end{equation}
	for each $x \in \mathbb{R}^{3N}$. In the notation above, $\mathcal{I}_{\left(l_1, l_2, \dots, l_{3N}\right)}(x)$ is the indicator function defined by
	\begin{equation*}
		\mathcal{I}_{\left(l_1, l_2, \dots, l_{3N}\right)}(x) = \begin{cases}
			1 & \text{if} \ \ \ x \in \prod_{i = 1}^{3N} \ \left[l_i h, (l_i+1)h  \right] \\
			0 & \text{else}.
		\end{cases}
	\end{equation*}
	Of course, we cannot sum~\eqref{eq:3N_histogram} over $l_i$ from negative infinity to positive infinity. Most summands will be zero anyway due the nature of the indicator function. One can truncate this expression by only summing over all bins that contain at least one particle in practice. The histogram evidently remains the same.
	
	\subsubsection{Calculating the Free Energy}
	With the histogram definition of $q_z$, we can estimate $\mu_\lambda(z)$ by~\eqref{eq:IS}. To calculate an estimate, we must draw $K$ i.i.d. samples $\{y_k\}_{k=1}^K$ from $H_z[\mathcal{X}_K](x)$. Suppose we have condensed the active bins to a list $\{B_j, K_j\}_{j=1}^J$, where $B_j$ is the actual bin set and $K_j$ is its particle count.
	
	We construct $y_k$ in two steps. First, draw a uniform number $u \sim U(0, 1)$ between $0$ and $1$. The particle $y_k$ falls in the $j$'th bin if
	\begin{equation*}
		\frac{\sum_{m=1}^{j-1} K_m}{K} \leq u <\frac{\sum_{m=1}^{j} K_m}{K}.
	\end{equation*}
	Second, draw a $3N$ dimensional sample from the uniform distribution on $B_j$. Carrying out this procedure for all samples $\{y_k\}_{k=1}^K$, the estimator for $N\mu_\lambda(z)$ reads
	\begin{equation}
		\tilde{\mu}_\lambda(z) = \frac{1}{K} \sum_{k=1}^K \frac{\mu_{\text{ext}}(z, y_k)}{H_z[\mathcal{X}_K](y_k)}.
	\end{equation}
	
	The histogram construction that we explained in this section is essentially inefficient. The reasons are twofold. In a practical implementation, we iterate over all reconstructed particles until we find a particle that has not been assigned a bin yet. Then, we loop over all remaining particles to calculate the bin count. This algorithm is essentially of $\mathcal{O}(K^2)$ complexity. There is no way to order the particles in advance to bring this complexity down.
	
	Secondly, drawing samples from the histogram is non-obvious because we generally do not know which bins are active. Constructing the list of active bins adds a layer of complexity to code. We present a different approach to building the histogram that mitigates these inefficiencies in the next section.
	
	\subsection{A Tensorized Histogram for Faster Estimates}  \label{subsec:tensor}
	There are many ways to define a histogram that can be used t estimate $\mu_\lambda$. The essential property is that the histogram must be explicitly normalized. In this section we construct an alternative tensorized histogram that mitigates the computational bottlenecks from previous section. This approach rests on the fact that if $H^{(i)}_z[\mathcal{X}^{(i)}_K](x^{(i)})$ is a histogram in dimension $i$, then
	\begin{equation} \label{eq:tensor_histogram}
		H_z[\mathcal{X}_K](x) = \prod_{i=1}^{3N} H^{(i)}_z[\mathcal{X}^{(i)}_K](x^{(i)}) 
	\end{equation}
	is a histogram in $3N$ dimensional space. We use the decomposition $x = \left(x^{(1)}, \dots, x^{(3N)}\right) \in \mathbb{R}^{3N}$.
	
	\subsubsection{Constructing the Histogram}
	
	Constructing a one dimensional histogram is then similar to section~\ref{subsec:histogram}, with the exception that we first sort the reconstructed particles $\mathcal{X}_K^{(i)}$ in ascending order. This brings the complexity of placing. particles in bins down to $\mathcal{O}(K  \log K)$. In dimension $i$, let $K_l^{(i)}$ be the number of  reconstructed particles $\mathcal{X}_K^{(i)}$ that fall into the interval $[l h, (l+1)h]$. Obviously, 
	\begin{equation*}
		\sum_{l=-\infty}^{\infty} K_l^{(i)} = K.
	\end{equation*}
	With these frequencies, the histogram in dimension $i$ reads
	\begin{equation} \label{eq:histogram1d}
		H_z^{(i)}[\mathcal{X}_K^{(i)}](x) = \frac{1}{K} \sum_{\ \in \mathbb{Z}} \mathcal{I}_{[lh, (l+1)h]}(x) K_l^{(i)}, \ x \in \mathbb{R}
	\end{equation}
	with $\mathcal{I}_{[lh, (l+1)h]}(x)$ the indicator function on the interval $[lh, (l+1)h]$. 
	
	\subsubsection{Calculating the Free Energy}
	We again draw $K$ i.i.d. samples $\{y_k\}_{k=1}^K$ from $H_z[\mathcal{X}_K]$. This is an easy task because of the tensorized definition~\eqref{eq:tensor_histogram}. Indeed, to generate a sample $y_k^{(i)}$ from the one-dimensional histogram $H_z^{(i)}[\mathcal{X}_K^{(i)}]$, draw a uniform number $u$ in $[0,1]$ and find the interval with number $m$ for which
	\begin{equation*}
		\frac{\sum_{l = -\infty}^{m-1} K_l^{(i)}}{K} \leq u <\frac{\sum_{l = -\infty}^{m} K_l^{(i)}}{K}.
	\end{equation*}
	Finally, draw a one dimensional microscopic sample $y_k^{(i)}$ uniformly in $[mh, (m+1)h]$. The combined sample in all $3N$ dimensions, $y_k = (y_k^{(1)}, \dots, y_k^{(3d)})$, is then distributed according to~\eqref{eq:tensor_histogram}. With these microscopic samples, an estimator for $\mu_{\lambda}$ finally reads
	\begin{equation} \label{eq:estimator_tensor}
		\tilde{\mu}_{\lambda}(z) = \frac{1}{K} \sum_{k=1}^K\frac{\mu_{\text{ext}}(z, y_k)}{H_z[\mathcal{X}_K](y_k)}.
	\end{equation}
	
	\begin{remark}
		When sampling from $H_z[\mathcal{X}_K]$ is on the critical computational path in an application, this step may be dropped entirely. It is possible to write an estimator for $\mu_{\lambda}$ using only the reconstructed samples $\{x_k\}_{k=1}^K$ as follows
		\begin{equation*}
			\tilde{\mu}_{\lambda}(z) = \frac{1}{K} \sum_{k=1}^K\frac{\mu_{\text{ext}}(z, x_k)}{H_z[\mathcal{X}_K](x_k)}.
		\end{equation*}
		That is, we use these microscopic samples to both construct ánd evaluate the histogram. This approach is not exact, but it can speed up computations if necessary.
	\end{remark}
	
	\subsubsection{Computing the Microscopic Acceptance Probability}
	Putting everything together from this section, our final pseudo-marginal approximation of the microscopic acceptance probability reads
	\begin{equation}
		\alpha_f(z_n, x_n, z', x')= \min \left\{1, \frac{\tilde{\mu}_\lambda(z') \bar{\mu}_0(z_n)}{\tilde{\mu}_\lambda(z_n) \bar{\mu}_0(z')} \right\}.
	\end{equation}
	We used the same notation for the microscopic and macroscopic particles as in section~\ref{sec:mMMCMC}.

	\subsection{Complete Algorithm} \label{subsec:algorithm}
\begin{figure}[h]
		\includegraphics[scale=0.6]{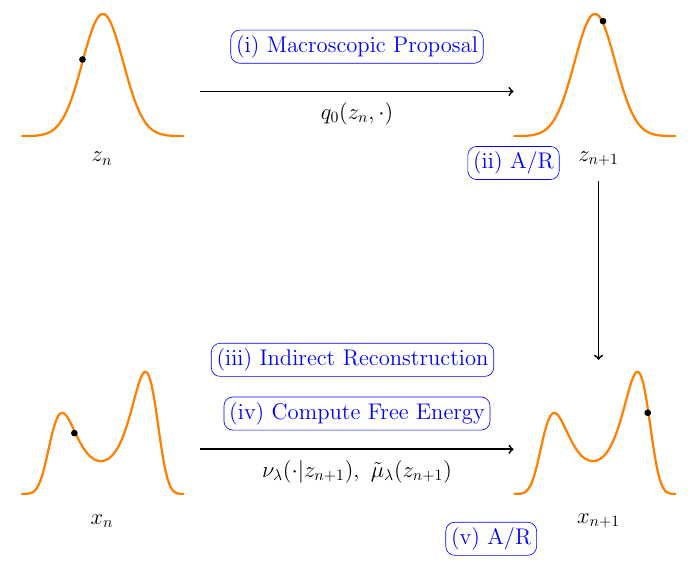}
		\caption{Illustration of the five steps in the mM-MCMC algorithm with pseudo-marginal approximation: macroscopic proposal, macroscopic accept/reject, indirect reconstruction, pseudo-marginal approximation of $\tilde{\mu}_{\lambda}$, and microscopic accept/reject.}
		\label{figuremMPM}
	\end{figure}

	The complete mM-MCMC method with pseudo-marginal approximation is illustrated in figure~\ref{figuremMPM}.
	
	\section{Numerical Results} \label{sec:numerics}
	Let us now review how the mM-MCMC method with pseudo-marginal approximation performs on molecular examples. We study the method on two molecules, butane and alkane. Butane is a simple molecule consisting of four chained carbon atoms, surrounded by hydrogen atoms. All chemical bonds are of covalent type. Alkane is a single chain that contains $N$ carbon atoms, and can be regarded as an extension of butane. The latter molecule allows us to study how the quality of pseudo-marginal approximation depends on the number of dimensions of the system. 
	
	To study the efficiency of the mM-MCMC method with pseudo-marginal approximation, we need a baseline method to compare possible efficiency gains. Our baseline will be the MALA method~\cite{xifara2014langevin,roberts1996exponential}. We explain the efficiency gain criterion further in section~\ref{subsec:efficiency_gain}. Afterwards, we carry out the experiments on butane and alkane in sections~\ref{subsec:butane} and~\ref{subsec:nalkane} respectively.
	
	\subsection{Efficiency Gain} \label{subsec:efficiency_gain}
	We first explain the Metropolis-Adjusted Langevin Algorithm (MALA) in some detail in section~\ref{subsubsec:mala}, after which we define the efficiency gain criterion in section~\ref{subsubsec:efficiency_gain}.
	
	\subsubsection{The Metropolis-Adjusted Langevin Algorithm (MALA)} \label{subsubsec:mala}
	MALA is a well known method to sample the time-invariant Gibbs distribution~\eqref{eq:gibbs}. Given the current microscopic sample $x_n \in \mathbb{R}^{3N}$, at time $t_n = n \ \delta t$, it generates the next sample as
	\begin{equation} \label{eq:mala}
		x_{n+1} = x_n - \delta t \nabla V(x_n) + \sqrt{2dt\beta} \eta_n,
	\end{equation}
	with $\eta_n$ distributed according to the $3N-$dimensional normal distribution with mean zero and unit covariance. We then accept $x_{n+1}$ with Metropolis-Hastings probability
	\begin{equation*}
		\alpha(x_n, x_{n+1}) = \min\left\{1, \frac{\mu(x_{n+1}) q(x_{n+1}, x_n)}{\mu(x_n) q(x_n, x_{n+1})}  \right\}.
	\end{equation*}
	In the above formula, $q(x_n, x_{n+1})$ is the transition probability distribution associated to~\eqref{eq:mala}.
	
	\subsubsection{Efficiency Gain Criterion} \label{subsubsec:efficiency_gain}
	Suppose we are interested in estimating the expected value of some functional $F: \mathbb{R}^{m} \to \mathbb{R}$ of reaction coordinate values with respect to the invariant measure
	\begin{equation} \label{eq:F}
		\mathbb{E}[F] = \int_{\mathbb{R}^n} F(\xi(x)) \ d\mu(x) = \int_{\mathbb{R}^n} F(z) \ d\mu_\lambda(z).
	\end{equation}
	By drawing random samples from $\mu$, through a Markov chain Monte Carlo method, we can obtain an estimate $\tilde{F}$ for $\mathbb{E}[F]$. To assess the accuracy of the MCMC method, one can perform $R$ independent runs, each with estimated value $\tilde{F}_i$, $i=1, \dots, R$ for~\eqref{eq:F}, and compute the Mean Squared Error (MSE).
	\begin{equation} \label{eq:mse_mMMCMC}
		\text{MSE}_{\text{mM-MCMC}} = \frac{1}{R} \sum_{i=1}^R \left( \tilde{F}_i -\mathbb{E}[F] \right)^2.
	\end{equation}
	In coming experiments, we will compare the performance of the mM-MCMC method to the microscopic MALA method. Therefore, we define the efficiency gain of mM-MCMC over MALA as
	\begin{equation} \label{eq:effgain}
		\frac{\text{MSE}_{\text{MALA}}}{\text{MSE}_{\text{mM-MCMC}}} \frac{T_{\text{MALA}}}{T_{\text{mM-MCMC}}}.
	\end{equation}
	Here, $\text{MSE}_{\text{MALA}}$ is defined analogously to~\eqref{eq:mse_mMMCMC}, and $T_{\text{MALA}}$ and $T_{\text{mM-MCMC}}$ are the respective CPU times of the MALA and mM-MCMC methods.
	
	In general, we can assume that the mean squared error of the mM-MCMC method will be lower than that of the MALA method. The reason is simple. The size of proposals by the MALA method are fundamentally limited by the fastest time scales. On the other hand, proposals by the mM-MCMC method are generated on the macroscopic level, bypassing this time-scale separation. As a result, the mM-MCMC method can explore the microscopic state space more efficiently. However, computing an approximation at every reconstruction step comes at a high cost. The MALA method will generally be a factor $2$ to $4$ times faster than the mM-MCMC for the same number of microscopic samples. Whether the mM-MCMC makes an efficiency gain over the MALA method comes down to a balance between these two factors.

	\subsection{Butane} \label{subsec:butane}
	The butane molecule is a single chain of four carbon atoms. The middle two atoms have two hydrogens and the other two have three. We are interested in the  torsion angle of the main chain, which will be our reaction coordinate. The potential energy surface (PES) for this molecule is
	\begin{equation} \label{eq:pes_butane}
		V(x)  = \sum_{CH_x - CH_y} \frac{1}{2} k_b (r(x) - r_0)^2 + \sum_{CH_x - CH_y - CH_z} \frac{1}{2} k_a (\theta(x) - \theta_0)^2 + A(\tau(x)).
	\end{equation}
	The first sum is taken over the three C-C bonds, the second over both C-C-C angles. The third term is the energy of the torsion angle. There are no electrostatic or Vanderwaals interactions~\cite{schappals2017round}. The constants have values $k_b / k_B  = 319225 K / \angstrom^2, \ k_a/k_B = 62500 \ K, r_0 = 1.540 \angstrom$ and $\theta_0 = 114^{\circ}$. Here, $k_B$ is the Boltzmann constant in appropriate units. The free energy of the reaction coordinate is
	\begin{equation} \label{eq:free_energy_butane}
		A(z) = A(\tau(x)) =  \sum_{i=0}^3 c_i \cos(\tau(x))^i,
	\end{equation}
	$c_0/k_B = 1031.36 \ K,\  c_1/k_B = 2037.82 \ K, \ c_2/k_B = 158.52 \ K, \ c_3/k_B = -3227.70 \ K.$ 
	
	We perform two experiments with mM-MCMC on butane. First, we show the histogram approximations by the mM-MCMC method with pseudo-marginal approximation to invariant probability distributions $\mu_\lambda$ and $\mu_0$. This experiment is in section~\ref{subsubsec:visual_inspection_butane}. Afterwards, we plot all estimated values for $Q_\lambda$ as a function of $z$, obtained during a single sample run of the mM-MCMC method in section~\ref{subsubsec:estimated_N_lambda}.

	\subsubsection{Inspection of mM-MCMC with Pseudo-Marginal Approximation} \label{subsubsec:visual_inspection_butane}
\begin{figure}[h]
		\centering
		\begin{subfigure}[b]{0.5\textwidth}
			\centering
			\includegraphics[width=\linewidth]{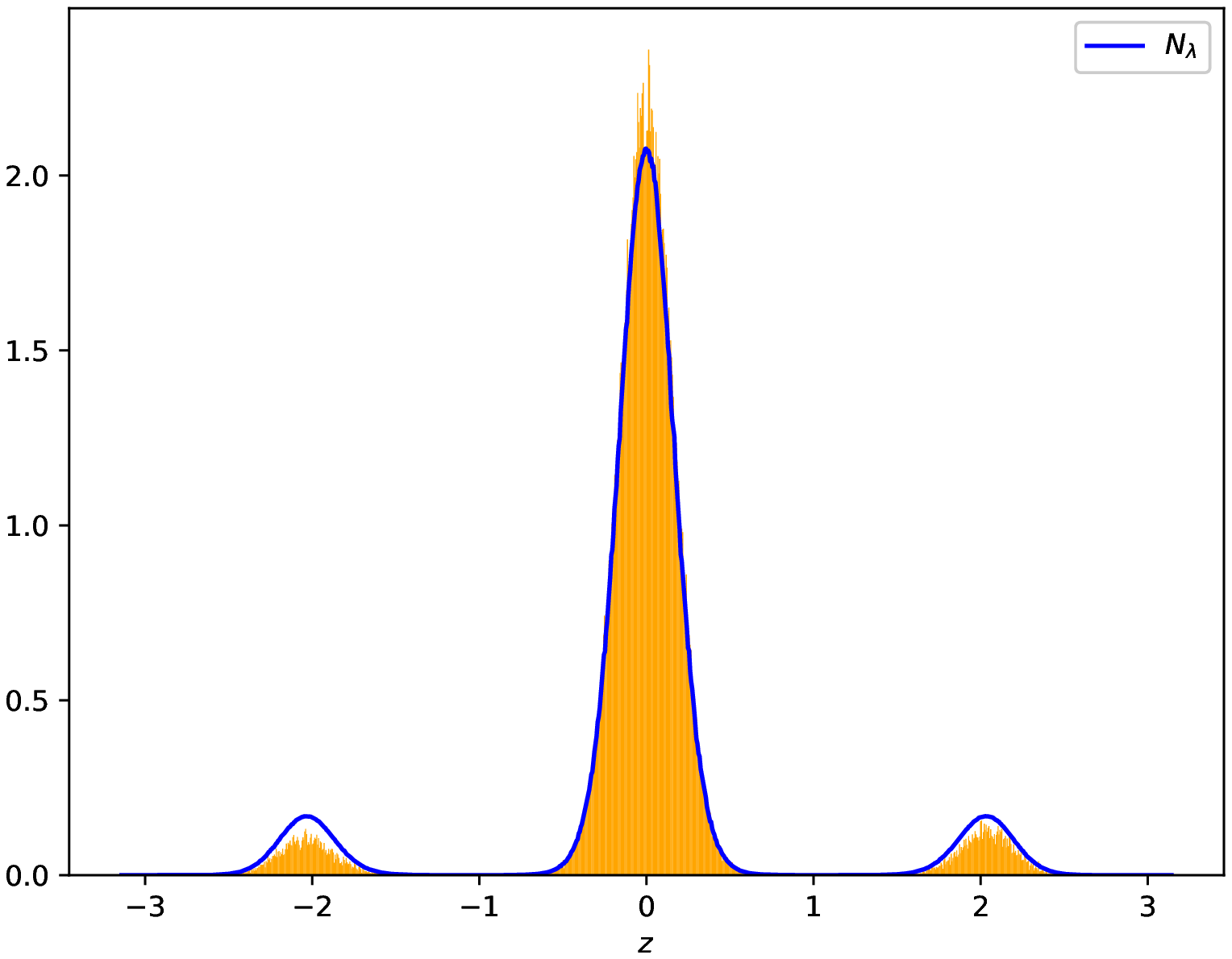}
		\end{subfigure}%
		\begin{subfigure}[b]{0.5\textwidth}
			\centering
			\includegraphics[width=\linewidth]{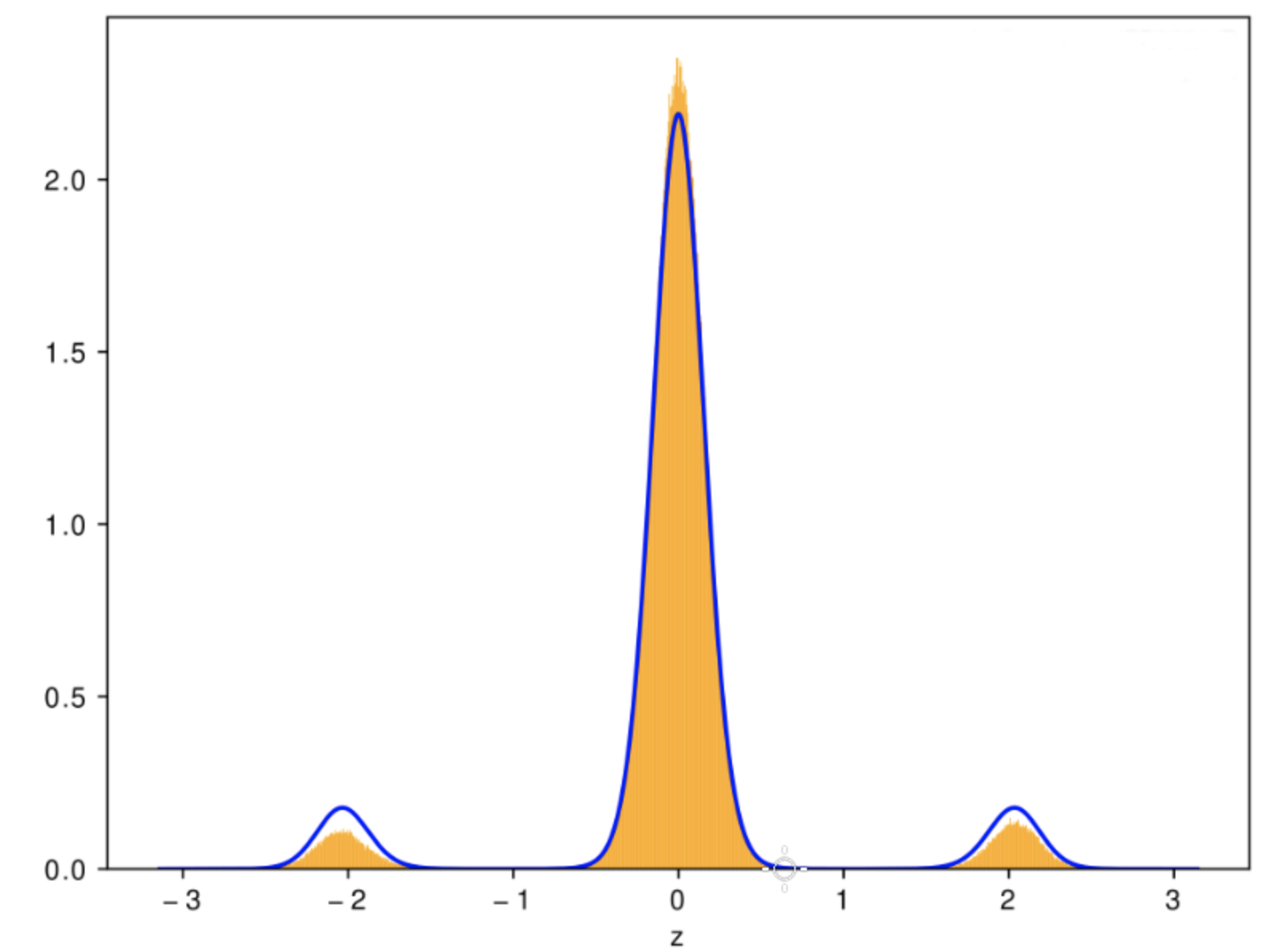}
		\end{subfigure}
		\caption{Histogram fit (orange) of a typical run of the mM-MCMC method with pseudo-marginal approximation on $\tilde{Q}_{\lambda}$ (blue, left) and the free energy (blue, right).}
		\label{fig:A_NL_z}
	\end{figure}

	We run the mM-MCMC algorithm with pseudo-marginal approximation for a total of $N=10^6$ sampling steps. On the macroscopic level, we use Brownian increments to sample reaction coordinate values with a time step $\Delta t = 0.001$. The biased simulation has parameters $\lambda = 2.0 \  k_b$, time step $\delta t = 0.01/\lambda$ and number  of steps $K=15$. The bin size of the histogram is $\sqrt{\frac{1}{2\lambda}}$. The torsion angles obtained by mM-MCMC with pseudo-marginal approximation are demonstrated in figure~\ref{fig:A_NL_z}. Left is the histogram of the reaction coordinate samples obtained on the macroscopic level, right those of the microscopic samples.
	
	Both histograms follow $\mu_\lambda(z)$ and respectively $\mu_0(z)$ well, even though there is an added variance because of the pseudo-marginal approximation.  The total macroscopic acceptance rate is $0.334$, and the microscopic acceptance rate is $0.75$. This number is smaller than the expected microscopic acceptance rate of $1.0$ due to the added variance by the pseudo-marginal estimation. Nevertheless, the approximation is unbiased.
	
	\subsubsection{Variance on Free Energy Estimates $\hat{Q}_\lambda(z)$ as a Function of Reaction Coordinate Value $z$} \label{subsubsec:estimated_N_lambda}
\begin{figure}[h]
		\centering
		\includegraphics[width=0.7\linewidth]{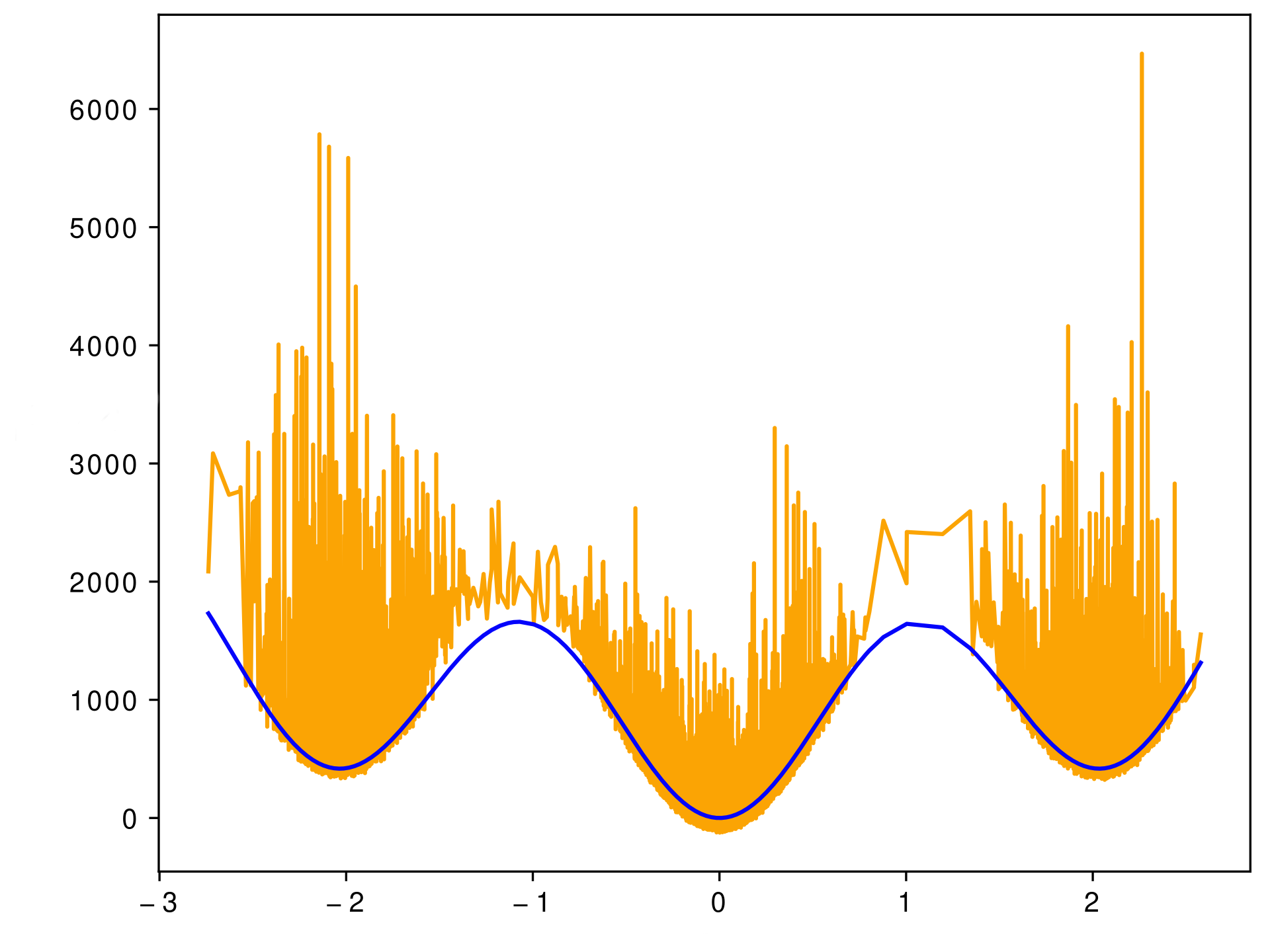}
		\caption{Value of $\tilde{Q}_\lambda(z)$ recorded during a mM-MCMC sample run (orange), as well as the exact value (blue).}
		\label{fig:nlambdapm}
	\end{figure}

	The induced variance on $\mu_\lambda$ decreases with $K$~\eqref{eq:IS}, but can vary with $z$. In this experiment, we run the mM-MCMC method with the same parameters as above, and record all values $\tilde{Q}_\lambda(z) $ seen during a sample run. We plot these values as well as the true free energy $Q_\lambda$ in figure~\ref{fig:nlambdapm}.
	
	There are two things to note. First, there is a very large variance on $\tilde{Q}_\lambda$. For two reaction coordinate values close together, the pseudo-marginal method computes $\mu_\lambda$ independently. These values can vary significantly, there is no continuity. Second, the actual approximate values follow the free energy profile well, there is no bias depending on $z$.
	
	\subsection{Alkane} \label{subsec:nalkane}
	An alkane molecule consists of a single chain of $N$ carbon atoms, with each carbon having two bonded hydrogen atoms. The carbon atoms at the end of the chain have an extra bonded hydrogen atom. It can be regarded as an extension of butane. Our motive for choosing this molecule is that it allows us to study how properties of the mM-MCMC method with pseudo-marginal approximation vary with dimensionality.
	
	With this objective in mind, we have designed four numerical experiments in this section. In the first experiment in section~\ref{subsubsec:avg_micro_acc}, we study the averaged microscopic acceptance rate as a function of $N$. We expect the microscopic acceptance rate to decrease as $N$ grows, but we do not know how fast. Then in section~\ref{subsubsec:acc_K_Ngroup}, we study the evolution of the averaged microscopic acceptance rate as function of the number of reconstructed particles $K$, for different values of $N$. Afterwards, we investigate the variance of the histogram estimator~\eqref{eq:estimator_tensor} for $\mu_\lambda$ as a function of $z$ and the number of reconstructed particles $K$. This experiment is in section~\ref{subsubsec:variance_K_Ngroup}. Finally, we combine everything by plotting the total efficiency gain on the expected mean and variance of $z$ as a function of $N$ and $K$ in section~\ref{subsubsec:MSE_K_Ngroup}. 
	
	\subsubsection{Averaged Microscopic Acceptance Rate as a Function of $N$} \label{subsubsec:avg_micro_acc}
\begin{figure}[h]
		\centering
		\includegraphics[width=0.7\linewidth]{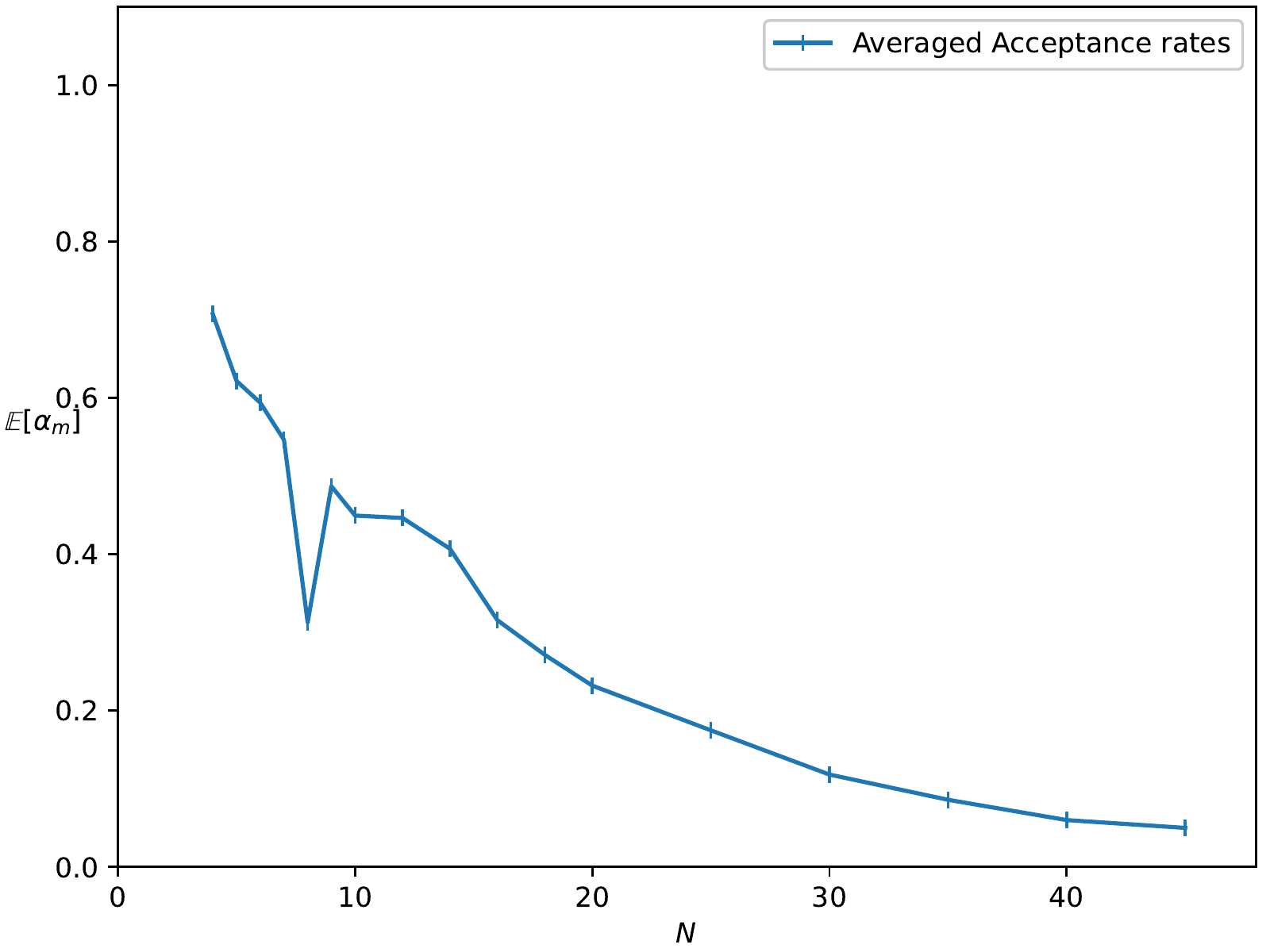}
		\caption{The averaged microscopic acceptance rate as a function of the number of carbon atoms (blue). Black lines indicate sampled values. The microscopic acceptance rate first decreases exponentially, then converges to an equilibrium value.}
		\label{fig:microacc}
	\end{figure}

	In this experiment, we analyse how the microscopic acceptance rate~\eqref{eq:alpha_pm} depends on $N$, the number of carbon atoms. The reaction coordinate is the torsion angle between the first four carbon atoms. For each value of $N$ between $4$ and $45$, we run the pseudo-marginal mM-MCMC algorithm for $10^5$ sampling steps. The macroscopic time step is $\Delta t = 0.001$, and number of indirect reconstruction steps $K=20$ with step size $\delta t = 0.01\ \lambda$. The constraining parameter is $\lambda = 2 k_b$. The average microscopic acceptance rate as a function of $N$ is pictured in figure~\ref{fig:microacc}.
	
	The acceptance rate decreases exponentially with $N$. This behaviour is typical for MCMC algorithms~\cite{cotter2013mcmc}, where the acceptance rate decreases exponentially with the number of dimensions of the system. Afterwards, however, the averaged microscopic acceptance rate converges to an equilibrium value. This value is small, but can be tweaked by the number of reconstruction steps. We explore this phenomenon further in the next experiment.
	
	\subsubsection{Averaged Microscopic Acceptance Rate as a Function of $K$} \label{subsubsec:acc_K_Ngroup}
\begin{figure}[h]
		\centering
		\includegraphics[width=0.7\linewidth]{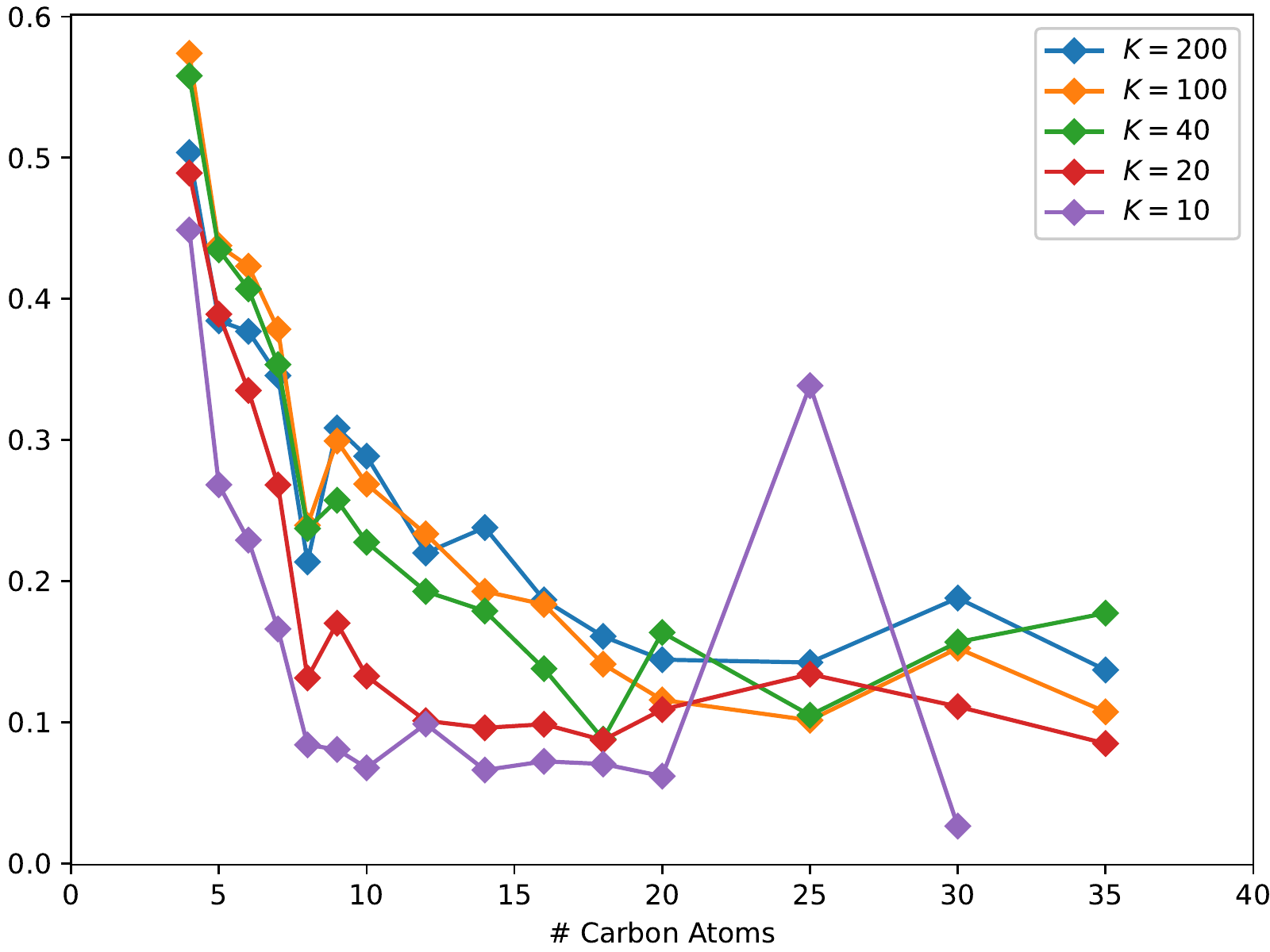}
		\caption{The averaged microscopic acceptance rate over $100$ independent runs as a function of $N$,. We show these curves for $K = 10$ (purple), $K=20$ (red). $K=40$ (green), $K=100$ (orange), and $K=200$ (blue).}
		\label{fig:acceptancekngroup}
	\end{figure}

	We saw in the previous experiment that the averaged microscopic acceptance rate first decreases, and then seemingly converges to some `limit' acceptance rate. The number of reconstructed particles was, however, on the low side. This analysis forms the basis of the current experiment. Our question is how the averaged microscopic acceptance rate depends on $K$, the number of reconstructed particles. The setup of this experiment is as follows. We run the mM-MCMC method with pseudo-marginal approximation for the same values of $N$ as above, and for $K = 10,20,40,100$ and $200$. The number of sampling steps is $10^5$. For each combination of parameters, we average the microscopic acceptance probability over $100$ independent runs. These curves are being shown on figure~\ref{fig:acceptancekngroup}.
	
	First of all, we again see the exponential decay of the averaged microscopic acceptance rate as $N$ gets larger. We already reached this conclusion in the last experiment. Second, the `limit' value increases with growing $K$. Indeed, the more reconstructed particles at each iteration of the mM-MCMC method, the lower the averaged variance on estimates of $\mu_\lambda$ will be. Then, according to~\eqref{eq:micro_acc_variance} the (averaged) microscopic acceptance probability should increase because it is fundamentally limited by uncertainty on $\mu_\lambda$.

	\subsubsection{Variance on $\mu_\lambda$ as a Function of $K$ and $z$} \label{subsubsec:variance_K_Ngroup}
	The analysis from the previous experiment relies on the fact that the variance on estimates of $\mu_\lambda(z)$ decreases with increasing $K$. We test this hypothesis here in this section. In this experiment, we plot the variance of statistical estimates of $\mu_\lambda$ as a function of $z$, and this for several values of $K$, the number of reconstructed particles. In a typical run of the mM-MCMC algorithm, there will never be two or more estimates at the same level of $z$ because the macroscopic proposals are everywhere between $-\pi$ and $\pi$. We need to aggregate estimates close to the current value of $z$, and weigh them according to their distance with respect to the current abscissa. A Gaussian weight function is a natural approach.

	Suppose $\{z_n, \tilde{Q}_\lambda(z_n)\}_{n=1}^N$ are the respective reaction coordinate values and estimated values for $Q_\lambda$ obtained during a run of the mM-MCMC method.  Then, for a fixed reaction coordinate value $z$, the mean and variance of $\tilde{Q}_\lambda$ are given by
	\begin{align} \label{eq:mean_variance_theory}
		\mathbb{E}  [ \tilde{Q}_\lambda | z ]&= \int_{-\pi}^{\pi} \int_{-\infty}^{\infty} \tilde{Q}_\lambda(u) \ d\tilde{Q}_\lambda(u) \ \delta_{u-z}(du) \\
		\mathbb{V} [\tilde{Q}_\lambda | z] &= \int_{-\pi}^{\pi}  \int_{-\infty}^{\infty} \left(\tilde{Q}_\lambda(u) - \mathbb{E} [\tilde{Q}_\lambda | z] \right)^2 \ d\tilde{Q}_\lambda(u) \ \delta_{u-z}(du)
	\end{align}
	As mentioned above, we then approximate the delta function by a sharp Gaussian factor centred at $z$ and with standard deviation $\varepsilon \ll 1$. The discrete, particle approximations to the integrals in~\eqref{eq:mean_variance_theory} then read
	\begin{align}
		m(z) &= \frac{1}{N} \sum_{n=1}^N \tilde{Q}_\lambda(z_n) \ \frac{1}{\sqrt{2\pi\varepsilon^2}}    \exp\left(-\frac{(z_n-z)^2}{2\varepsilon^2}\right) \\
		\sigma^2(z)  &=  \frac{1}{N} \sum_{n=1}^N \left(\tilde{Q}_\lambda(z_n) - m(z) \right)^2 \ \frac{1}{\sqrt{2\pi\varepsilon^2}}    \exp\left(-\frac{(z_n-z)^2}{2\varepsilon^2}\right).
	\end{align}
	
\begin{figure}[h]
		\centering
		\includegraphics[width=0.7\linewidth]{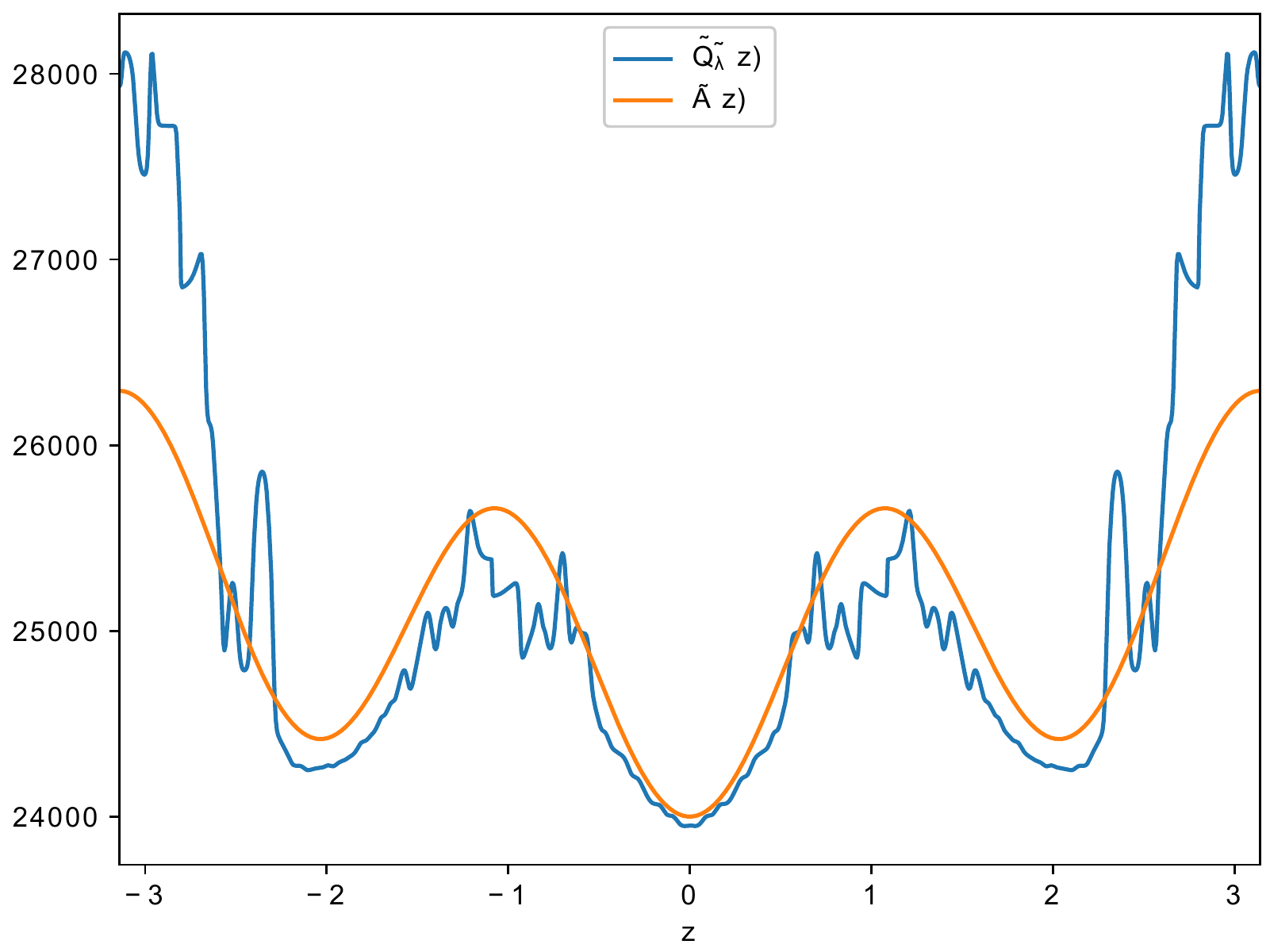}
		\caption{Average of estimates for $\tilde{Q}_\lambda(z)$ averaged over $100$ independent runs (blue), and the free energy $A(z)$ (orange).}
		\label{fig:expectationnlambdaz}
	\end{figure}

	In the remainder of this section, we show two figures. The first figure depicts the aggregated mean $m(z)$ for $z$ between $-\pi$ and $\pi$, averaged over $100$ independent runs of the mM-MCMC algorithm with $10^5$ samples. The other parameters are the same as in section~\ref{subsubsec:visual_inspection_butane}. The number of carbon atoms is $8$. Further parameters of the mM-MCMC method are $K = 20$, $\Delta t = 0.001$,  $\lambda = 2 k_b$, and $\delta t = 0.01  / \lambda$. The standard deviation of the Gaussian weight function is $\varepsilon=0.02$. This plot, together with the free energy~\eqref{eq:free_energy_butane} are shown on figure~\ref{fig:expectationnlambdaz}.
	
\begin{figure}[h]
		\centering
		\includegraphics[width=0.8\linewidth]{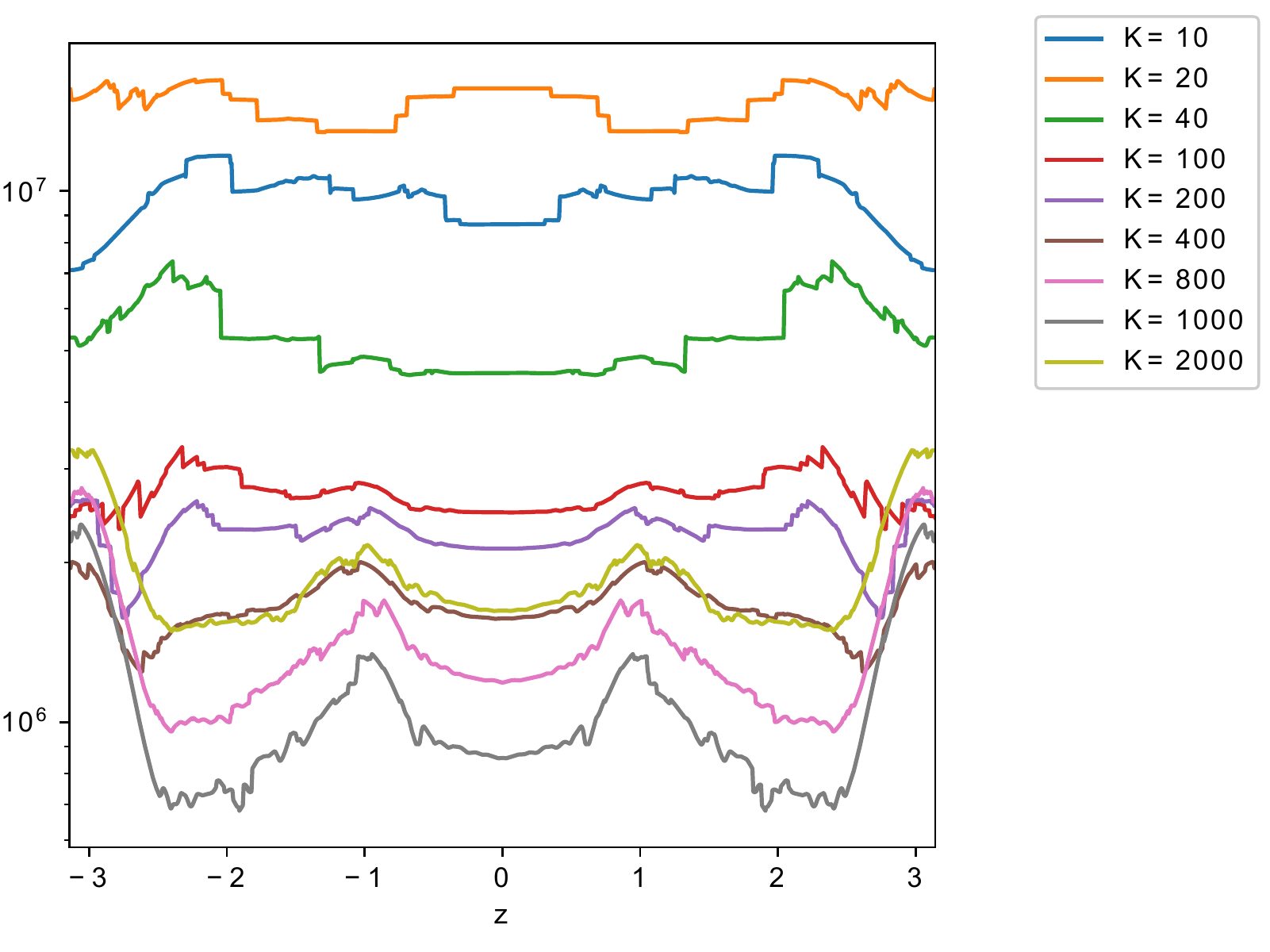}
		\caption{Variance of estimates for $\tilde{Q}_\lambda(z)$ averaged over $100$ independent runs for several values of $K$, shown in the legend.}
		\label{fig:variancenlambdak}
	\end{figure}

	In the second figure, we show the variance on estimates $\tilde{\mu}_\lambda(z_n)$ for different values of $K$ These values are $10, 20, 40, 100, 200, 400, 1000$, and $2000$. All other parameters remain the same. These variance plots are shown in figure~\ref{fig:variancenlambdak}.

	The results are clear. On figure~\ref{fig:expectationnlambdaz} the averaged curve $m(z)$ follows $A(z) \approx Q_\lambda(z)$ closely, especially in regions where the density of reaction coordinate values are large. These regions are centred at $-2, 0$ and $2$. Only at the ends of the figure do the curves diverge because there are very few samples there to get an accurate estimate.
	
	We see a similar pattern on figure~\ref{fig:variancenlambdak}. First of all, the variance on estimates $\tilde{Q}_\lambda$ does indeed decrease when $K$ increases. This decrease happens at a steady pace, note the vertical axis uses a log-scale. This decrease confirms our statement from section~\ref{subsubsec:acc_K_Ngroup}. Second, the variance is lowest in regions of high density, and the variance is largest at both ends of the figure, especially for large $K$.

	\subsubsection{Efficiency gain as a Function of $K$ and $N$} \label{subsubsec:MSE_K_Ngroup}
	The only factor we have not taken into account yet is the cost of reconstruction. We have shown in section~\ref{subsec:tensor} that computing $\mu_\lambda$ has a log-linear complexity in $K$. Although the total sampling accuracy grows with larger $K$, so does the total cost. It is unclear how the sampling accuracy grows with the number of reconstruction steps. We test this dependency here.
	
\begin{figure}[h]
		\centering
		\includegraphics[width=0.7\linewidth]{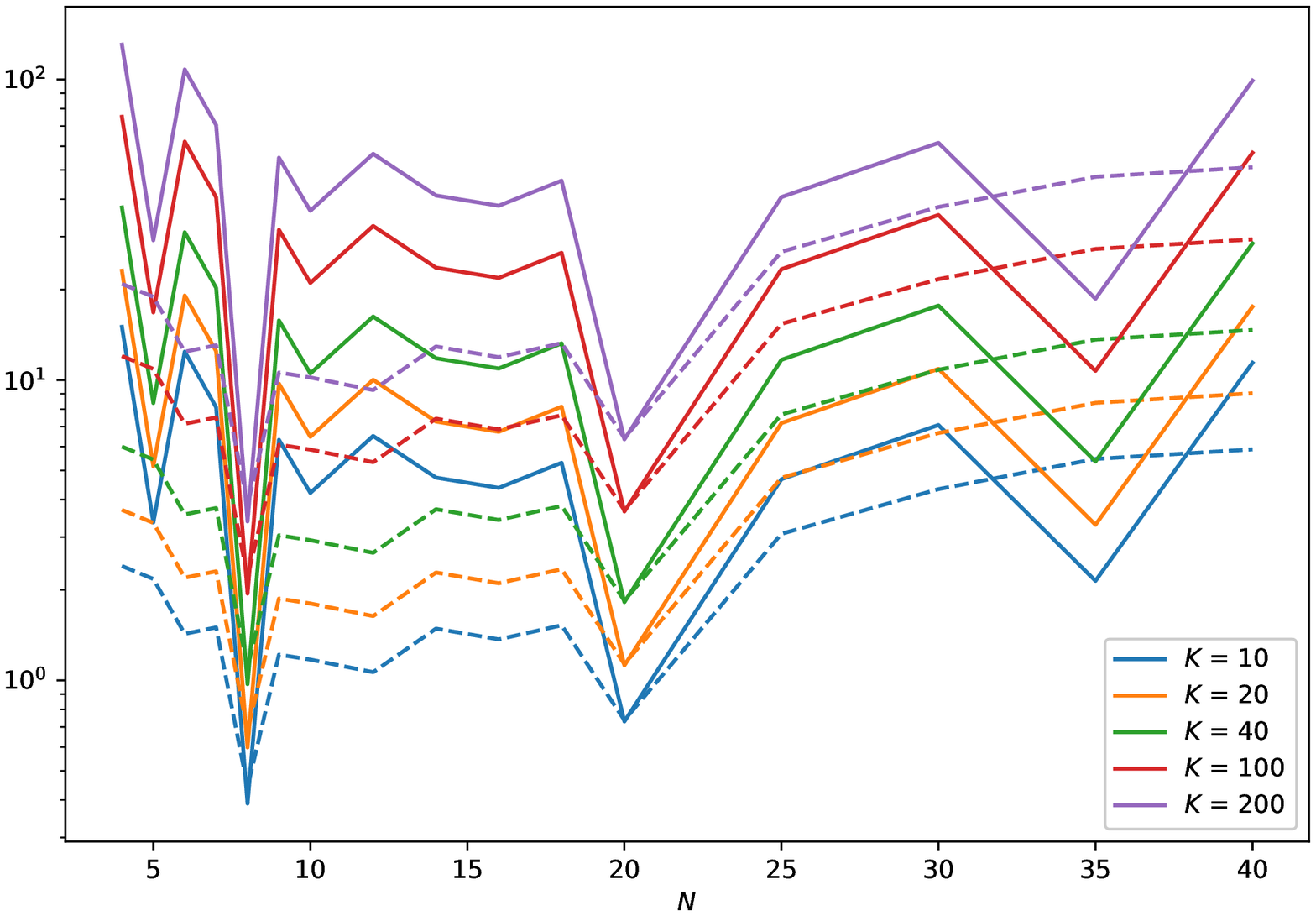}
		\caption{Efficiency gain of the mM-MCMC method with pseudo-marginal approximation over the microscopic MALA method on the mean (full) and variance (dashed) of the torsion angle for several values of $K$.}
		\label{fig:efficiencygain}
	\end{figure}

	In this experiment, we plot the total efficiency gain of the mM-MCMC method with pseudo-marginal approximation over the microscopic MALA method~\ref{subsubsec:mala} for several values of $N$ and $K$. We specifically measure the efficiency gain with respect to the estimated mean and variance of the reaction coordinate. In the language of~\eqref{eq:F}, the functions $F$ are respectively
	\begin{equation*}
		\begin{cases}
			F_{\text{mean}}(z) &= z \\
			F_{\text{variance}}(z) &= z^2.	
		\end{cases}
	\end{equation*}
	The setup for this experiment is as follows. For each value of $N$, we run the mM-MCMC method for $10^5$ sampling steps. The macroscopic time step is $\Delta t = 0,001$, and $\lambda = 2k_b$. The number of reconstruction steps $K$ are $10, 20, 40, 100$ and $200$.  The number of sampling steps of the microscopic MALA method is also $10^5$, and its time step is $\delta t = 0.01/\lambda$. We compute the efficiency gain over $100$ independent runs. The results are shown on figure~\ref{fig:efficiencygain}.
	
	The results are unexpected. The mM-MCMC method is on average $\mathbb{E}[\alpha_{CG}] K \log K$ times slower per microscopic sample than the MALA method. For $K=20$ and an averaged macroscopic acceptance rate of $0.3$, it comes down to a factor $7.8$ slower per sample. Still, the efficiency gain on $F_{\text{mean}}$ and $F_{\text{variance}}$ increase with increasing $K$, a very positive result. The reason can only be that mM-MCMC gets accurate fast with increasing $K$, enough to offset the ratio of runtimes in~\eqref{eq:effgain}. We hope the reader will be convinced of the benefit of the mM-MCMC method.

	\section{Conclusion} \label{sec:conclusion}
	We developed a pseudo-marginal approach to estimate the free energy inside the micro-macro Markov chain Monte Carlo method. The free energy of the reaction coordinate appears inside the microscopic acceptance probability, and is hard and computationally expensive to estimate. The pseudo-marginal method makes a statistical estimate to the free energy on the fly by means of importance sampling. This estimate is unbiased and uses few microscopic samples.  We implemented pseudo-marginal mM-MCMC on two molecules, butane and alkane. From the experiments, we saw that the pseudo-marginal methods computes accurate approximations to the free energy with only a small amount of samples. Furthermore, the variance on these free energy estimates decreases steadily for larger amounts of samples. When we studied the whole mM-MCMC method with pseudo-marginal approximation, we noticed the averaged microscopic acceptance rate increases for larger numbers reconstructed samples, a result consistent with the theory. In the final experiment, we computed the overall efficiency gain on the mean and variance of the free energy by the mM-MCMC method over the microscopic MALA method. The efficiency gain unexpectedly increases with larger numbers of reconstructed, even though the cost per microscopic sample increases as well.

	\section*{Acknowledgements}
	This work was funded by the Flemish Fund for Scientific Research (FWO) with Grant 1179820N for fundamental research. The resources and services used in this work were provided by the VSC (Flemish Supercomputer Center), funded by the Research Foundation - Flanders (FWO) and the Flemish Government.
	\bibliographystyle{plain}
	\bibliography{bib}
	
\end{document}